\documentclass{amsart}
\usepackage{amsmath}
\usepackage{amssymb}
\usepackage{amscd}
\usepackage{latexsym}
\usepackage{graphics, epsfig}
\usepackage{curves}

\theoremstyle{plain}
\newtheorem{thm}{Theorem}[section]
\newtheorem{lem}[thm]{Lemma}

\newtheorem{prop}[thm]{Proposition}

\theoremstyle{definition}

\def \CPb {\overline{\mathbf{CP}}^{2}}
\def \CP {{\mathbf{CP}}^{2}} 
\def \R {\mathbf{R}}
\def \Z {\mathbf{Z}}
\def \Sig{\Sigma}

\def \g {\gamma}
\def \d {\delta}

\def \lam {\lambda}
\def \G {\Gamma}

\def \o {\omega}

\def \bd {\partial}
\def \x {\times}
\def \- {\setminus}
\def \C {\subset}
\def \ve {\varepsilon}

\def \ssw {\text{SW}}

\def \DD {\Delta}

\def\hk{\widehat{k}}

\begin{document}

\baselineskip.5cm

\title [Double node neighborhoods] {Double node neighborhoods and families of simply connected 4-manifolds with ${\bf{b^+=1}}$}
\author[Ronald Fintushel]{Ronald Fintushel}
\address{Department of Mathematics, Michigan State University \newline
\hspace*{.375in}East Lansing, Michigan 48824}
\email{\rm{ronfint@math.msu.edu}}
\thanks{The first author was partially supported NSF Grant DMS0305818
and the second author by NSF Grant DMS0204041}
\author[Ronald J. Stern]{Ronald J. Stern}
\address{Department of Mathematics, University of California \newline
\hspace*{.375in}Irvine,  California 92697}
\email{\rm{rstern@math.uci.edu}}


\maketitle

\section{Introduction\label{Intro}}

A basic question of $4$-manifold topology is whether the complex projective plane, $\CP$, admits exotic smooth structures. Thus one is interested in knowing the smallest $m$ for which $\CP\# m\,\CPb$ admits an exotic smooth structure. In the late 1980's, Dieter Kotschick   \cite{K} proved that the Barlow surface, which was known to be homeomorphic to $\CP\#\, 8\CPb$, is not diffeomorphic to it. In following years the subject of simply connected smooth $4$-manifolds with $b^+=1$  languished because of a lack of suitable examples.  However, largely due to a beautiful paper of Jongil Park \cite{P}, who found the first examples of exotic simply connected $4$-manifolds with $b^+=1$ and $b^-=7$, the past year has found renewed interest in this subject. Peter Ozsvath and Zoltan Szabo proved that Park's manifold is minimal \cite{OS}, and Andras Stipsicz and Szabo used a technique similar to Park's to construct an exotic manifold with $b^+=1$ and $b^-=6$ \cite{SS}.

The rational elliptic surface $E(1)\cong \CP\# 9\,\CPb$ admits infinitely many distinct smooth structures. Up to now it has not been known whether $\CP\# m\,\CPb$ can have an infinite family of smooth structures  when $m<9$. It is the purpose of this paper to introduce a new technique which we use to show that for $m=6$, $7$, $8$, $\CP\# m\,\CPb$ does have an infinite family of smooth structures. The construction of these examples to some extent follows the ideas of Park \cite{P} after performing suitable knot surgeries \cite{KL4M} on $E(1)$. The essential new idea is that after certain knot surgeries on $E(1)$, one is able to find an immersed sphere of self-intersection $-1$ representing the `pseudo-section'. This is accomplished by studying `double node neighborhoods' as described below.

Shortly after the appearance of a preliminary version of this article, Park, Stipsicz, and Szabo \cite{PSS} used the techniques described herein to give even better examples: an infinite family of pairwise nondiffeomorphic smooth $4$-manifolds all homeomorphic to $\CP\# 5\,\CPb$. In the final section of this paper we give a quick construction of such a family of examples. 

\section{Seiberg-Witten invariants, knot surgery, and rational blowdowns}

\subsection{Seiberg-Witten invariants} Let $X$ be a simply connected oriented 
4-mani-fold with $b_X^+=1$ with a given orientation of $H^2_+(X;\R)$ and a given metric $g$. The Seiberg-Witten invariant depends on the metric $g$ and a self-dual 2-form as follows. 
There is a unique $g$-self-dual harmonic 2-form $\o_g\in H^2_+(X;\R)$ with $\o_g^2=1$ and corresponding to the positive orientation. Fix a characteristic homology class $k\in H_2(X;\Z)$.  Given a pair $(A,\psi)$, where
$A$ is a connection in the complex line bundle whose first Chern class is the Poincar\'e dual $\hk=\frac{i}{2\pi}[F_A]$ of $k$ and $\psi$ a section of the bundle $W^+$ of self-dual spinors for the associated $spin^{\,c}$ structure, the perturbed Seiberg-Witten equations are:
\begin{gather*} 
D_A\psi = 0 \\
F_A^+  = q(\psi)+i\eta \notag\label{SWeqn}
\end{gather*}
where $F_A^+$ is the self-dual part of the curvature $F_A$,
$D_A$ is the twisted Dirac operator, $\eta$ is a
self-dual 2-form on $X$, and
$q$ is a quadratic function. Write $\ssw_{X,g,\eta}(k)$ for the
corresponding invariant. As the pair
$(g,\eta)$ varies, $\ssw_{X,g,\eta}(k)$ can change only at those pairs
$(g,\eta)$ for which there are solutions with $\psi=0$. These 
solutions occur for pairs $(g,\eta)$ satisfying $(2\pi\hk+\eta)\cdot\o_g=0$.
This last equation defines a wall in $H^2(X;\R)$. 

The point $\o_g$ determines a component of the double cone consisting of elements of $H^2(X;\R)$ of positive square. We prefer to work with $H_2(X;\R)$. The dual component is determined by the Poincar\'e dual $H$ of $\omega_g$.  (An element $H'\in H_2(X;\R)$ of positive square lies in the same component as $H$ if $H'\cdot H>0$.) If
$(2\pi \hk+\eta)\cdot\o_g\ne 0$ for a generic $\eta$, $\,\ssw_{X,g,\eta}(k)$ is
well-defined, and its value depends only on the sign of $(2\pi \hk+\eta)\cdot\o_g$. Write $\ssw_{X,H}^+(k)$ for $\ssw_{X,g,\eta}(k)$ if 
$(2\pi \hk+\eta)\cdot\o_g>0$ and $\ssw_{X,H}^-(k)$ in the other case.

The invariant $\ssw_{X,H}(k)$ is defined by $\ssw_{X,H}(k) =\ssw_{X,H}^+(k)$ if 
$(2\pi \hk)\cdot\o_g>0$, or dually, if $k\cdot H>0$, and $\ssw_{X,H}(k) =\ssw_{X,H}^-(k)$ if $H\cdot k <0$. The wall-crossing formula \cite{KM,LL} states that if $H', H''$ are elements of positive square in $H_2(X;\R)$ with $H'\cdot H>0$ and $H''\cdot H>0$, then if $k\cdot H' <0$ and $k\cdot H''>0$,
\[ \ssw_{X,H''}(k) - \ssw_{X,H'}(k) = (-1)^{1+\frac12 d(k)}\]
where $d(k)=\frac14(k^2-(3\,\text{sign}+2\,\text{e})(X))$ is the formal dimension of the Seiberg-Witten moduli spaces.

Furthermore, in case $b^-\le 9$, the wall-crossing formula, together with the fact that $\ssw_{X,H}(k)=0$ if $d(k)<0$, implies that $\ssw_{X,H}(k) = \ssw_{X,H'}(k)$ for any $H'$ of positive square in $H_2(X;\R)$ with $H\cdot H'>0$. So in case $b^+_X=1$ and $b^-_X\le 9$, there is a well-defined Seiberg-Witten invariant, $\ssw_X(k)$.

\subsection{Rational blowdowns} Let $C_{p}$ be the smooth $4$-manifold obtained by plumbing $(p-1)$ disk bundles over the $2$-sphere according to
the diagram

\begin{picture}(100,60)(-90,-25)
 \put(-12,3){\makebox(200,20)[bl]{$-(p+2)$ \hspace{6pt}
                                  $-2$ \hspace{96pt} $-2$}}
 \put(4,-25){\makebox(200,20)[tl]{$u_{0}$ \hspace{25pt}
                                  $u_{1}$ \hspace{86pt} $u_{p-2}$}}
  \multiput(10,0)(40,0){2}{\line(1,0){40}}
  \multiput(10,0)(40,0){2}{\circle*{3}}
  \multiput(100,0)(5,0){4}{\makebox(0,0){$\cdots$}}
  \put(125,0){\line(1,0){40}}
  \put(165,0){\circle*{3}}
\end{picture}

\noindent Then the classes of the $0$-sections have self-intersections $u_0^2=-(p+2)$ and $u_i^2=-2$, $i=1,\dots,p-2$. The boundary of $C_p$ is the lens space 
$L(p^2, 1-p)$ which bounds a rational ball $B_p$ with $\pi_1(B_p)=\Z_p$ and $\pi_1(\bd B_p)\to \pi_1(B_p)$ surjective. If $C_p$ is embedded in a $4$-manifold $X$ then the rational blowdown manifold $X_{(p)}$ of  \cite{rat} is obtained by replacing $C_p$ with $B_p$, i.e., $X_{(p)} = (X\- C_p) \cup B_p$. 

The rationally blown down manifold $X_{(p)}$ shares many properties with $X$. For example, if $X$ and $X\- C_p$ are simply connected, then so is $X_{(p)}$. Also, the Seiberg-Witten invariants of $X$ and $X_{(p)}$ can be compared. The homology of $X_{(p)}$ can be identified with the orthogonal complement of the classes $u_i$, $i=0,\dots,p-2$ in $H_2(X;\Z)$, and then each characteristic element $k\in H_2(X_{(p)};\Z)$ has a lift $\widetilde{k} \in H_2(X;\Z)$ which is characteristic and for which the dimensions of moduli spaces agree, $d_{X_{(p)}}(k)=d_X(\widetilde{k})$. It is proved in \cite{rat} that if $b^+_X>1$ then $\ssw_{X_{(p)}}(k)=\ssw_X(\widetilde{k})$. In case $b^+_X=1$, if $H\in H_2^+(X;\R)$ is orthogonal to all the $u_i$ then it also can be viewed as an element of $H_2^+(X_{(p)};\R)$, and 
$\ssw_{X_{(p)},H}(k)=\ssw_{X,H}(\widetilde{k})$.

\subsection{Knot surgery} Let $X$ be a $4$-manifold which contains a homologically essential torus $T$ of self-intersection $0$, and let $K$ be a knot in $S^3$. Let $N(K)$ be a tubular neighborhood of $K$ in $S^3$, and let $T\x D^2$ be a tubular neighborhood of $T$ in $X$. Then the knot surgery manifold $X_K$ is defined by 
\[ X_K = (X\- (T\x D^2))\cup (S^1\x (S^3\- N(K))\]
The two pieces are glued together in such a way that 
the homology class $[{\text{pt}}\x \bd D^2]$ is identified with $[{\text{pt}}\x \lam]$ where $\lam$ is the class of a longitude of $K$. This latter
condition does not, in general, completely determine the diffeomorphism type of $X_K$; however if we take $X_K$ to be any manifold constructed in this fashion and if, for example, $T$ has a cusp neighborhood, then the Seiberg-Witten invariant of $X_K$ is completely determined by the Seiberg-Witten invariant of $X$ and the Alexander polynomial of $K$
\cite{KL4M}. Furthermore, if $X$ and $X\- T$ are simply connected, then so is $X_K$.

\section{Double node neighborhoods}

A simply connected elliptic surface is fibered over $S^2$ with smooth fiber a torus and with singular fibers. The most generic type of singular fiber is a nodal fiber, which is an immersed  $2$-sphere with one transverse positive double point. A nearby smooth fiber contains a vanishing cycle. This vanishing cycle is a nonseparating loop on the smooth fiber and the nodal fiber is obtained by collapsing this vanishing cycle to a point to create a transverse self-intersection. The vanishing cycle  bounds a `vanishing disk', a disk of relative self-intersection $-1$ with respect to the framing of its boundary given by pushing the loop off itself on the smooth fiber. The monodromy of the singular fiber is the diffeomorphism of the smooth fiber which describes the torus bundle over a small circle in the base $S^2$ bounding a disk whose only singular point is the image of the nodal fiber. In this case the monodromy is a Dehn twist around the vanishing cycle.

Consider $E(1)$ which admits an elliptic fibration with 12 nodal fibers, six of which have monodromy (a Dehn twist around) $a$ and six of which have monodromy (a Dehn twist around) $b$ where $a$ and $b$ represent a standard basis of $H_1(T^2;\Z)$. (Note that in the mapping class group of $T^2$, $ab$ has order $6$.) A standard technique for finding spheres of self-intersection $-2$ in a elliptic fibration is to find vanishing cycles on smooth fibers that vanish in two different nodal fibers. In the example of $E(1)$, the vanishing cycle $a$ on a smooth fiber bounds vanishing disks which are determined by a choice of node with monodromy $a$. Two such choices give a sphere of self-intersection $-2$, the union of the two vanishing disks.

A {\em{double node neighborhood}} $D$ is a fibered neighborhood of an elliptic fibration which contains exactly two nodal fibers with the same monodromy. If $F$ is a smooth fiber of $D$, there is a vanishing class $a$ that bounds vanishing disks in the two different nodes, and these give rise to a sphere $V$ of self-intersection $-2$ in $D$. 

We now investigate the result of a knot surgery along a regular fiber in a double node neighborhood. Let $K$ be any genus one knot in $S^3$ which contains a nonseparating loop $\G$ on a minimal genus Seifert surface $\Sig$ such that $\G$ satisfies the following two properties:
\begin{itemize}
\item[(i)] $\G$ bounds a disk in $S^3$ which intersects $K$ in exactly two points.
\item[(ii)] The linking number in $S^3$ of $\G$ with its pushoff on $\Sig$ is $+1$.
\end{itemize}   
It follows from these properties that $\G$ bounds a punctured torus in $S^3\- K$. Examples of $(K,\G)$ which satisfy (i) and (ii) are the twist knots $T(n)$ with Alexander polynomials $\DD_{T(n)}(t) = nt-(2n-1)+nt^{-1}$: 

\centerline{\unitlength 1cm
\begin{picture}(6.75,4)
\put (2.75,2.5){\oval(1.5,1)[l]}
\put (2.75,2.5){\oval(2.5,2)[l]}
\put (3.75,2.5){\oval(1.5,1)[r]}
\put (3.75,2.5){\oval(2.5,2)[r]}
\put (2.75,3.49){\line(1,0){.25}}
\put (2.75,2.99){\line(1,0){.25}}
\put (3.5,3.49){\line(1,0){.25}}
\put (2.95,3.175){\oval(1,.65)[rt]}
\put (3.6,3.325){\oval(1,.65)[lb]}
\put (3.5,2.99){\line(1,0){.25}}
\put (2.75,1.3){\framebox(1,1)}
\put (2.85,2){\Small{$2n-1$}}
\put (2.85,1.7){\Small{RH $\frac12$-}}
\put (2.85,1.35){\Small{twists}}
\put (1.9,.9) {\small{$T(n) =$ twist knot}}
\end{picture}}
\vspace{-.3in}\noindent where $\G$ is the loop which runs through both half-twists in the clasp. 

Now consider the result of knot surgery in $D$ using a knot $K$ with a loop $\G$ on its Seifert surface so that $(K,\G)$ satisfies (i) and (ii). In the knot surgery construction, one has a certain freedom in choosing the gluing of $S^1\x (S^3\- N(K))$ to  $D\- N(F)$. We are free to make any choice as long as a longitude of $K$ is sent to the boundary circle of a normal disk to $F$. We choose the gluing so that the class of a meridian $m$ of $K$ is sent to the class of $a \x \{pt\}$ in $H_1(\bd(D\- N(F));\Z)=H_1(F\x \bd D^2;\Z)$. 

Before knot surgery, the fibration of $D$ has a section which is a disk. The result of knot surgery is to remove a smaller disk in this section and to replace it with the Seifert surface of $K$. Call the resulting relative homology class in $H_2(D_K,\bd;\Z)$ the {\it{pseudo-section}}. Thus in $D_K$, the result of knot surgery, there is a punctured torus $\Sig_s$ representing the pseudo-section. The loop $\G$ sits on $\Sig_s$ and by (i) it bounds a twice-punctured disk $\DD$ in $\{pt\}\x \bd(S^3\- N(K))$ where $\bd\DD=\G\cup m_1\cup m_2$ where the $m_i$ are meridians of $K$. The meridians $m_i$ bound disjoint vanishing disks $\DD_i$ in $D\- N(F)$ since they are identified with disjoint loops each of which represents the class of $a \x \{pt\}$ in $H_1(\bd(D\- N(F));\Z)$. Hence in $D_K$ the loop $\G\C \Sig_s$ bounds a disk $U=\DD\cup \DD_1\cup \DD_2$. By construction, the relative self-intersection of $U$ relative to the framing given by the pushoff of $\G$ in $\Sig_s$ is $+1-1-1=-1$. (This uses (ii).)

\begin{lem} \label{imm} Let $(K,\G)$ be a knot together with a loop $\G$ on its Seifert surface which satisfies (i) and (ii), and let $D$ be a double node neighborhood. Then in $D_K$ the pseudo-section is represented by an immersed disk $\Lambda$ with one positive double point.
\end{lem}
\begin{proof} Since $\G$ is nonseparating in $\Sig_s$, surgery on it kills $\pi_1(\Sig_s)$. This surgery may be performed in $D_K$ by removing an annular neighborhood of $\G$ and replacing it with a pair of disks $U'$, $U''$ as obtained above. This is precisely the complex-algebraic model of a nodal intersection. So the resultant surface is as claimed.
\end{proof}

\section {Families of $4$-manifolds} \label{constr}

In this section we will construct a family of mutually nondiffeomorphic simply connected $4$-manifolds with $b^+=1$ and $b^-=6$. For the reader who has read Park's lovely paper \cite{P}, our construction closely parallels the construction there; however we first do a knot surgery on $E(1)$ and then need to blow up one less time to get started. Begin with $E(1)$. As Park points out, $E(1)$ has an elliptic fibration with 5 singular fibers, an $\widetilde{E}_6$ fiber and 4 nodal fibers. The nodal fibers can be chosen to consist of a pair of fibers with monodromy $a$ and two others with monodromies $b$ and $a^{-1}ba$ where $a$ and $b$ generate $\pi_1(T^2)$. (We obtain an elliptic fibration of $E(1)$ with this collection of singular fibers by factoring the monodromy $(ab)^6=(ab)^4a^2(a^{-1}ba)b$, and noting that $(ab)^4$ is the monodromy of an $\widetilde{E}_6$ fiber.)

Thus we can find a double node neighborhood $D\C E(1)$ containing the two nodal fibers with monodromy $a$ and so that $E(1)\- D$ contains an $\widetilde{E}_6$ fiber and the two remaining nodal fibers $F_1$ and $F_2$. 

Let $\eta$ be the class of a line in $\CP$ and $\ve_i$, $i=1,\dots,9$, be the classes of the exceptional curves in the elliptic surface $E(1) = \CP\# 9\,\CPb$. Then each $\ve_i$ is a section of the elliptic fibration, and the fiber $F$ represents $3\eta-\ve_1-\dots-\ve_9$.

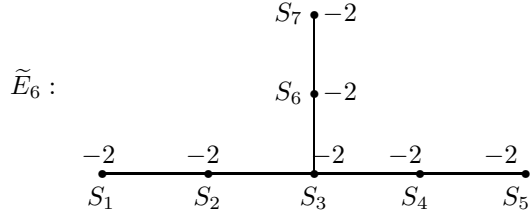
\centerline{
 \begin{picture}(300,90)(-65,-10)
   \put(2,3){\makebox(200,20)[bl]{$-2$ \hspace{20pt} $-2$
             \hspace{28pt} $-2$ \hspace{13pt}$-2$\hspace{20pt} $-2$}}
   \put(4,-25){\makebox(200,20)[tl]{$S_{1}$ \hspace{23pt} $S_{2}$
               \hspace{23pt} $S_{3}$ \hspace{21pt} $S_{4}$
               \hspace{21pt} $S_{5}$}}
   \put(75,58){$S_{7}$  \hspace{1pt} $-2$}
   \put(75,28){$S_{6}$  \hspace{1pt} $-2$}
   \multiput(90,60)(0,-30){2}{\line(0,-1){30}}
   \multiput(90,60)(0,-30){2}{\circle*{3}}
   \multiput(10,0)(40,0){4}{\line(1,0){40}}
   \multiput(10,0)(40,0){5}{\circle*{3}}
\put(-25,30){$\widetilde{E}_6:$}
 \end{picture}
 \label{E6}}

\vspace*{.1in}

Let $K$ be the twist knot $K=T(n)$, and let $Y_n = E(1)_K$ be the result of knot surgery using $K$ and identifying a smooth fiber $F$ of the elliptic fibration with $T=S^1\x m$ in $S^1\x (S^3\- N(K))$ (and where $m$ denotes a meridian to $K$ in $S^3$). We can do the surgery inside the double node neighborhood $D$; so $Y_n= (E(1)\- D)\cup D_K$. 

It follows from Lemma~\ref{imm} that $Y_n$ contains an immersed sphere $S$ of self-intersection $-1$ which is obtained from replacing a disk in the  exceptional curve representing $\ve_9$ with the immersed disk $\Lambda$ given by Lemma~\ref{imm}. In $E(1)\- D = Y_n\- D_K$ the two nodal fibers $F_1$ and $F_2$ each intersect $S$ transversely in a single positive point. Now blow up $Y_n$ three times, at the double points of $S$, $F_1$ and $F_2$, and let $Z_n$ be the blowup, $Z_n = Y_n\# 3\,\CPb$, with exceptional curves $E_i$, $i=0,1,2$.

In $Z_n$ we have the configuration consisting of the embedded sphere $S'$ of self-intersection $-5$ representing $S-2E_0$ and the embedded self-intersection $-4$ spheres $F_i'$ representing $F-2E_i$, $i=1,2$. These three spheres meet transversely in a pair of positive double points which can be locally smoothed to obtain an embedded sphere $R$ of self-intersection $-9$ and which represents the class \[ u_0 = S+F_1+F_2-2(E_0+E_1+E_2)= S+2T-2(E_0+E_1+E_2) \] 

As Park points out in \cite{P}, the classes of the spheres $S_i$ in $\widetilde{E}_6$ are: 
\begin{align*}  [S_1] &=  \ve_4 - \ve_7\\   [S_2] &=  \ve_1 - \ve_4\\
 [S_3] &=\eta -\ve_1 - \ve_2 -\ve_3  \\
 [S_4] &= \ve_2 - \ve_5\\ [S_5] &= \ve_5 - \ve_9\\ 
 [S_6] &= \ve_3 - \ve_6\\  [S_7] &= \ve_6 - \ve_8
 \end{align*}
Furthermore, $R$ intersects $S_5$ in a single positive point; so we obtain the configuration $C_7\C Z_n$ whose homology classes are: 
\begin{eqnarray*}
&u_0&\!\!\!\!=[R]=[S+2T-2(E_0+E_1+E_2)]\\
&u_1&\!\!\!\!=[S_5]= \ve_5 - \ve_9\\
&u_2&\!\!\!\!=[S_4]= \ve_2 - \ve_5\\
&u_3&\!\!\!\!=[S_3]= \eta -\ve_1 - \ve_2 -\ve_3\\
&u_4&\!\!\!\!=[S_2]=  \ve_1 - \ve_4\\
&u_5&\!\!\!\!=[S_1]=  \ve_4 - \ve_7
\end{eqnarray*}
Recall that while in $Z_n$ the notation $\ve_i$ is meaningless, the classes which we have listed, such as $\ve_5 - \ve_9$ make perfect sense, since they are represented by embedded spheres in $Y_n\- D_K = E(1)\- D$. 

Let $X_n = (Z_n\- C_7)\cup B_7$ denote the manifold obtained by rationally blowing down the configuration $C_7$. The boundary of $C_7$ is the lens space $L(49,-6)$ and $\pi_1(L(49,-6))=\Z_{49}$ maps onto $\pi_1(B_7)=\Z_7$.
A normal circle to $S_3$ is a loop representing a generator of $\pi_1(\bd C_7)$, and this loop is the boundary of the disk $S_6\- C_7$. Hence $X_n$ is simply connected. It is clear that $b^+_{X_n}=1$ and $b^-_{X_n}=6$. Thus

\begin{prop} For each $n$, the manifold $X_n$ is homeomorphic to $\CP\#6\,\CPb$.\qed
\end{prop}

\section{The Seiberg-Witten invariants of $X_n$}\label{sw}

In this section we shall compute the Seiberg-Witten invariants of the manifolds $X_n$, $n>0$, and see that they are mutually nondiffeomorphic. We first recall the result of \cite{S} and \cite{KL4M} concerning the manifolds $Y_n$, $n>0$. Recall that since $b^-_{Y_n}=9$, the manifold $Y_n$ has a well-defined Seiberg-Witten invariant. 

\begin{lem}[\cite{S,KL4M}] The Seiberg-Witten invariants of $Y_n$, $n>0$, are given by: $|\ssw_{Y_n}(\pm T)| = n$ and for every other class $L$, \ $\ssw_{Y_n}(L)=0$. \qed
\end{lem}

In $E(1)$ the sphere $S_{\eta}$ representing $\eta$ intersects the fiber $F$ in 3 points. After knot surgery, $\eta$ gives rise to a class $h$ in $Y_n$ of genus 3 that has $h^2=1$ and $h\cdot T =3$. (The three normal disks to $F$ that lie in $S_{\eta}$ are replaced by genus one Seifert surfaces of $K=T(n)$.) Since the Seiberg-Witten invariant of $Y_n$ is well-defined, $\ssw_{Y_n,h}(L)=\ssw_{Y_n}(L)$ for all characteristic $L\in H_2(Y_n;\Z)$.

The homology $H_2(Y_n;\R)$ embeds naturally in $H_2(Z_n;\R)$, and the blowup formula of \cite{Turkey} implies that in $Z_n$, the only classes $L$ for which the Seiberg-Witten invariants $\ssw_{Z_n,h}(L)$ can be nonzero are 
$L=\pm T-\d_1E_1-\d_2E_2-\d_3E_3$ where the $\d_i$ are odd integers (since $L$ is characteristic). If $d(L)$ is the dimension of the moduli space of solutions to the Seiberg-Witten equations corresponding to $L$ (and $h$), then $4d(L)=L^2-(3\,\text{sign}(Z_n)+2\,e(Z_n)) = -(\d_1^2+\d_2^2+\d_3^2)+3$. Hence the blowup formula, and the fact that $\ssw_{Z_n,h}(L)=0$ if $d(L)<0$,  implies:

\begin{lem} For $n>0$,  $|\ssw_{Z_n,h}(\pm T\pm E_0\pm E_1\pm E_2)| = n$,   and $\ssw_{Z_n,h}(L)=0$ for all other classes $L$. \qed
\end{lem}

The wall-crossing formula then implies:

\begin{prop} \label{ssw} For any class $H\in H_2(Z_n;\R)$ satisfying $H\cdot h>0$ and $H^2>0$, and for any characteristic $L\in H_2(Z_n;\R)$ with $L^2\ge -3$, if $L=\pm T\pm E_1\pm E_2\pm E_3$ we have
\[ \ssw_{Z_n,H}(L)=\begin{cases} \pm n\ \ &{\text{if the signs of $H\cdot L$ and $h\cdot L$ agree}}\\
\pm n\pm 1 &  {\text{if the signs of $H\cdot L$ and $h\cdot L$ do not agree}}
\end{cases}\]
and if $L\ne \pm T\pm E_1\pm E_2\pm E_3$ we have
\[ \ssw_{Z_n,H}(L)=\begin{cases} 0\ \ &{\text{if the signs of $H\cdot L$ and $h\cdot L$ agree}}\\
\pm 1 &  {\text{if the signs of $H\cdot L$ and $h\cdot L$ do not agree}} \qed
\end{cases}\]
\end{prop}

We are interested in the Seiberg-Witten invariants $\ssw_{X_n}(K)$. The homology $H_2(X_n;\Z)$ of $X_n$ is naturally identified with the orthogonal complement to $\{u_i|i=0,\dots,5\}$ in $H_2(Z_n;\Z)$. It is proved in \cite{rat} that for each
characteristic element $K\in H_2(X_n;\Z)$ there is a lift $\widetilde{K}$ which is a characteristic element of $H_2(Z_n;\Z)$ satisfying 
\begin{multline*} d_{X_n}(K)=\frac14(K^2-(3\,\text{sign}(X_n)+2\,e(X_n))) =\\
\frac14(\widetilde{K}^2-(3\,\text{sign}(Z_n)+2\,e(Z_n)))=d_{Z_n}(\widetilde{K}) 
\end{multline*}
The lift $\widetilde{K}$ is obtained by extending the restriction of $K$ to $\bd B_7=\bd C_7$ over $C_7$ as a characteristic vector whose self-intersection with respect to the relative intersection form on $C_7$ is $-6$. It is then showed in \cite{rat} that
\[ \ssw_{X_n}(K) = \ssw_{Z_n,H}(\widetilde{K}) \]
where $H\in H_2(Z_n;\R)$ is any class satisfying $H\cdot h>0$, $H^2>0$, and  $H\cdot u_i=0$ for $i=0,\dots,5$.

Such a class $H$ is given by 
\[ H= 7h - 2(\sum_{i=1}^9e_i ) - e_3-E_0-E_1-E_2\]
Analogously to the situation for $h$, the classes $e_i \in H_2(Y_n;\Z)\C H_2(Z_n;\Z)$ are obtained from the exceptional curves in $E(1)$ after knot surgery. Each is represented by a torus of self-intersection $-1$ and the intersection properties of $h, e_1,\dots,e_9$ in $Y_n$ are the same as those of $\eta,\ve_1,\dots,\ve_n$ in $E(1)$ (and note $\ve_9=[S]$). It follows that $H\cdot h= 7$ and $H^2 =5$. Furthermore, $H\cdot u_i=0$ for $i=0,\dots,5$. So the class $H$ may be used to determine the Seiberg-Witten invariants of $X_n$. 

The intersection matrix of $C_7$ in terms of the basis of $H_2(C_7,\bd;\Z)$ given by the dual classes $\g_i$ to the $u_i$ is the inverse of the plumbing matrix of $C_7$ (see \cite{rat}). Using this, one easily checks that of the classes 
$T\pm E_0\pm E_1\pm E_2$, only $\widetilde{K}_0=T +E_0+ E_1 + E_2$ restricts to $C_7$ with self-intersection equal to $-6$. In fact, $\widetilde{K}_0$ restricts to $H_2(C_7,\bd;\Z)$ as $7\g_0$. The induced class $K_0 \in H_2(X_n;\Z)$ is characteristic, and $\ssw_{X_n}(K_0)=\ssw_{Z_n,H}(T+E_0+E_1+E_2)$. Calculating intersections, we have 
$H\cdot (T+E_0+E_1+E_2)= 5$ and $h\cdot (T+E_0+E_1+E_2)=3$. It then follows from Proposition~\ref{ssw} that 
\[ |\ssw_{X_n}(\pm K_0)|=|\ssw_{Z_n,H}(T+E_0+E_1+E_2)| = |\ssw_{Z_n,h}(T+E_0+E_1+E_2)|= n.\]

\begin{prop} \label{key} For $n\ge 2$, there is (up to sign) a unique class $K_0\in H_2(X_n;\Z)$ whose Seiberg-Witten invariant has absolute value unequal to $0$ or $1$, and we have $|\ssw_{X_n}(\pm K_0)|= n.$ \qed
\end{prop}

As a consequence, we have the family of smooth manifolds promised in the introduction:

\begin{thm} The manifolds $X_n$ ($n>0$) are all homeomorphic to $\CP\# 6\,\CPb$, and no two of these manifolds are diffeomorphic.
Furthermore, for $n\ge 2$, the $X_n$ are all minimal.
 \end{thm}
\begin{proof} Any diffeomorphism $X_n\to X_m$ must take $K_0$ to a class in $H_2(X_m;\Z)$ whose Seiberg-Witten invariant has absolute value $n$. It then follows from Proposition~\ref{key} that $n=m$. Since $\CP\# 6\,\CPb$ admits a metric of positive scalar curvature, its Seiberg-Witten invariant is identically $0$. Hence it is not diffeomorphic to any $X_n$, $n>0$.

According to the blowup formula \cite{Turkey}, if $X_n$ were diffeomorphic to $X'\#\CPb$, the basic classes of $X_n$ would come in pairs $k\pm E$ where $E$ is the exceptional class coming from the blowup, $k\in H_2(X';\Z)$, and 
$\ssw_{X_n}(k\pm E)=\ssw_{X'}(k)$. Then, if $n\ge 2$,  we would have $K_0 =k\pm E$ and $-K_0 = k\mp E$ for some $k$ since these are the only two classes in $H_2(X_n;\Z)$ whose Seiberg-Witten invariants are not $0$ or $\pm 1$. Thus $(2 K_0)^2=
(K_0-(-K_0))^2 = (2E)^2=-4$. But $K_0^2= 3$; so this is not possible.
\end{proof}

Note that the manifolds $X_n$ have no underlying symplectic structure for $n \ge 2$. They provide the first examples of nonsymplectic $4$-manifolds homeomorphic to $\CP\# m\,\CPb$ with $m < 9$.

\section{Families with $b^-=7$ and $8$}

One can easily modify the construction of  Section \ref{constr} to obtain infinite families of simply connected minimal $4$-manifolds with $b^+=1$ and $b^-=7$ or $8$. In the first case simply blow up twice rather than three times, using only one nodal fiber to construct a sphere of self-intersection $-7$ in $Y_n\# 2\,\CPb$, and add on enough spheres in $\widetilde{E}_6$ to form the configuration $C_5$. To obtain families with $b^+=8$, only blow up the double point on the pseudo-section in $Y_n$ to get a sphere of self-intersection $-5$ in $Y_n\#\CPb$, then add a sphere to get the configuration $C_3$. Calculations of Seiberg-Witten invariants are analogous to those of Section \ref{sw}.

\section{Families with $b^-=5$} 

Park, Stipsicz, and Szabo \cite{PSS} have recently shown how to utilize the techniques of this paper to produce infinite families of mutually nondiffeomorphic $4$-manifolds, all homeomorphic to $\CP\#5\,\CPb$. In this section, we outline one such family. This construction was also discovered by Park, Stipsicz, and Szabo. The construction begins by noting (as in \cite{PSS}) that the standard elliptic fibration on $E(1)$ has monodromy $(ab)^6$, where $a$ and $b$, as in Section~4, correspond to Dehn twists around standard generators of $\pi_1(T^2)$. Using the braid relation, $aba=bab$, one sees that $(a^3b)^3$ defines an elliptic fibration on $E(1)$, and then so does $a^6(a^{-3}ba^3)(bab^{-1})^2b^2(b^{-1}ab)$. This means that there is an elliptic fibration on $E(1)$ whose singular fibers are an $I_6$ fiber (whose monodromy is $a^6$ \cite{BPV}), and $6$ nodal fibers, whose monodromies are 
$a^{-3}ba^3$, two with $bab^{-1}$, two with $b$, and one with $b^{-1}ab$. 

Since we have two pairs of nodal fibers with parallel vanishing cycles, we can form a double node neighborhood containing each pair. Perform a knot surgery in each double node neighborhood; say using the twist knot $T(1)$ in one neighborhood and $T(n)$ in the other. In each neighborhood we need to be careful to perform the knot surgery so that the meridian of the knot is identified with the vanishing cycle. Let $V_n$ be the resultant manifold. Note that $V_n$ is simply connected: For if $F_1$, $F_2$ are the fibers on which we perform knot surgery, then $\pi_1(E(1)\-  (F_1\cup F_2))$ is normally generated by normal circles to the removed fibers, and these are identified with longitudes of $T(1)$ and $T(n)$. Thus $\pi_1(V_n)$ is normally generated by the images of $\pi_1(S^1\x (S^3\- T(i)))$, for $i=1, n$.
These, in turn are normally generated by the images of $S^1\x \{\text{pt}\}$ and the meridians. However, because of the singular fibers, $\pi_1({\text{fiber}})\to \pi_1(E(1)\- (F_1\cup F_2))$ is the trivial map. So the classes in question die in $\pi_1(E(1)\- (F_1\cup F_2))$, and we see that $\pi_1(V_n)=1$. 
Using \cite{KL4M}, one computes the Seiberg-Witten invariants of $V_n$: Up to sign the only classes in $H_2(V_n;\Z)$ with nontrivial Seiberg-Witten invariants are $T$ and $3T$, and $|\ssw_{V_n}(\pm 3T)|=n$, $|\ssw_{V_n}(\pm T)|=2n-1$.

To construct the examples of manifolds with $b^-=5$, use the two double node neighborhoods to get a representative of the pseudo-section of $V_n$ which is an immersed sphere with two positive double points. If we blow up these two points we get an embedded sphere $S$ of self-intersection $-9$ in $W_n=V_n\# 2\,\CPb$. By removing one of the spheres in the $I_6$ configuration, we get $5$ spheres which, along with $S$ form the configuration $C_7$. As in Section \ref{sw}, we can see that the only classes $L\in H_2(W_n;\Z)$ which have $\ssw_{W_n,h}(L)\ne 0$ and whose restriction to $C_7$ has self-intersection 
$-6$ are $L=\pm(3T+E_0+E_1)$. If we let $Q_n$ denote the result of rationally blowing down $C_7$ in $W_n$, then an argument exactly as in Section \ref{sw} shows:  

\begin{thm}[{\it{cf.}} Park, Stipsicz, and Szabo \cite{PSS}]  The manifolds $Q_n$ ($n>0$) are all homeomorphic to $\CP\# 5\,\CPb$, and no two of these manifolds are diffeomorphic.
Furthermore, for $n\ge 2$, the $Q_n$ are all minimal.
\end{thm}
\noindent (Note that to see that $Q_n$ is simply connected, we use either of the two nodal fibers which do not live in the double node neighborhoods.)


\begin{thebibliography}{9999}

\bibitem[BPV]{BPV} W. Barth, C. Peters, and A. Van de Ven, ``Compact Complex Surfaces,'' Springer-Verlag, 1984.

\bibitem[FS1]{rat} R. Fintushel and R. Stern, {\em Rational blowdowns of smooth $4$-manifolds}, Jour. Diff. Geom. {\bf 46} (1997), 181--235.

\bibitem[FS2]{Turkey} R. Fintushel and R. Stern, {\em Immersed spheres in
4-manifolds and the immersed Thom conjecture}, Turkish J. Math. {\bf 19}
(1995), 145--157.

\bibitem[FS3]{KL4M} R. Fintushel and R. Stern, {\em Knots, links, and
$4$-manifolds}, Invent. Math. {\bf 134} (1998), 363--400.

\bibitem[K]{K}  D. Kotschick, {\em On manifolds homeomorphic to
${\mathbf CP}^2 \sharp 8{\overline{{\mathbf CP}}^2}$},
Invent. Math. {\bf 95} (1989), 591--600.

\bibitem[KM]{KM} P. Kronheimer and T. Mrowka, {\em The genus of embedded
surfaces in the projective plane}, Math. Research Letters {\bf 1} (1994),
797--808.

\bibitem[LL]{LL} T.J. Li and A. Liu, {\em Symplectic structure on ruled surfaces and a generalized adjunction formula}, Math. Research Letters {\bf 2} (1995),
453--471.

\bibitem[OS]{OS} P. Ozsvath and Z. Szabo, {\em On Park's exotic smooth four-manifolds}, preprint, http://front.math.ucdavis.edu/math.GT/0411218.
 
\bibitem[P]{P} J. Park, {\em Simply connected symplectic $4$-manifolds with $b_2^+=1$ and $c_1^2=2$}, Invent. Math. (to appear), http://front.math.ucdavis.edu/math.GT/0311395.

\bibitem[PSS]{PSS} J. Park, A. Stipsicz and Z. Szabo, {\em Exotic smooth structures on $\CP\#5\,\CPb$}, preprint, http://front.math.ucdavis.edu/math.GT/0412216.

\bibitem[SS]{SS} A. Stipsicz and Z. Szabo, {\em An exotic smooth structure on $\CP\# 6\,\CPb$}, preprint, http://front.math.ucdavis.edu/math.GT/0411258.

\bibitem[S]{S} Z. Szabo, {\em Exotic $4$-manifolds with $b^+_2=1$}, Math. Res. Letters {\bf 3} (1996), 731--741.


\end{thebibliography}
\end{document}